\documentclass[11pt,a4paper]{amsart}
\usepackage{graphics,array,amsmath,amssymb,epsf,graphicx}
\usepackage[export]{adjustbox}
\usepackage{esvect}

\newtheorem{thm}{Theorem}
\newtheorem{lem}[thm]{Lemma}
\newtheorem{prop}[thm]{Proposition}
\newtheorem{cor}[thm]{Corollary}

\renewcommand{\qed}{\hfill$\Box$}

\begin{document}
\title{Yamada Polynomial and associated link of $\theta$-curves}
\author{Youngsik Huh}
\address{Department of Mathematics,
College of Natural Sciences, Hanyang University, Seoul 04763,
Korea} \email{yshuh@hanyang.ac.kr}




\begin{abstract}
The discovery of polynomial invariants of knots and links, ignited by V. F. R. Jones, leads to the formulation of polynomial invariants of spatial graphs. The Yamada polynomial, one of such invariants, is frequently utilized for practical distinguishment of spatial graphs. Especially for $\theta$-curves, the polynomial is an ambient isotopy invariant after a normalization. On the other hand, to each $\theta$-curve, a 3-component link can be associated as an ambient isotopy invariant. The benefit of associated links is that invariants of links can be utilized as invariants of $\theta$-curves.

In this paper we investigate the relation between
the normalized Yamada polynomial of $\theta$-curves and the Jones polynomial of their associated links, and show that the two polynomials are equivalent for brunnian $\theta$-curves as a corollary. For our purpose the Jaeger polynomial of spatial graphs is observed, a specialization of which is equivalent to the Yamada polynomial.
\end{abstract}

\maketitle



\section{Introduction}
\begin{figure}
\centering
\includegraphics[width=6cm]{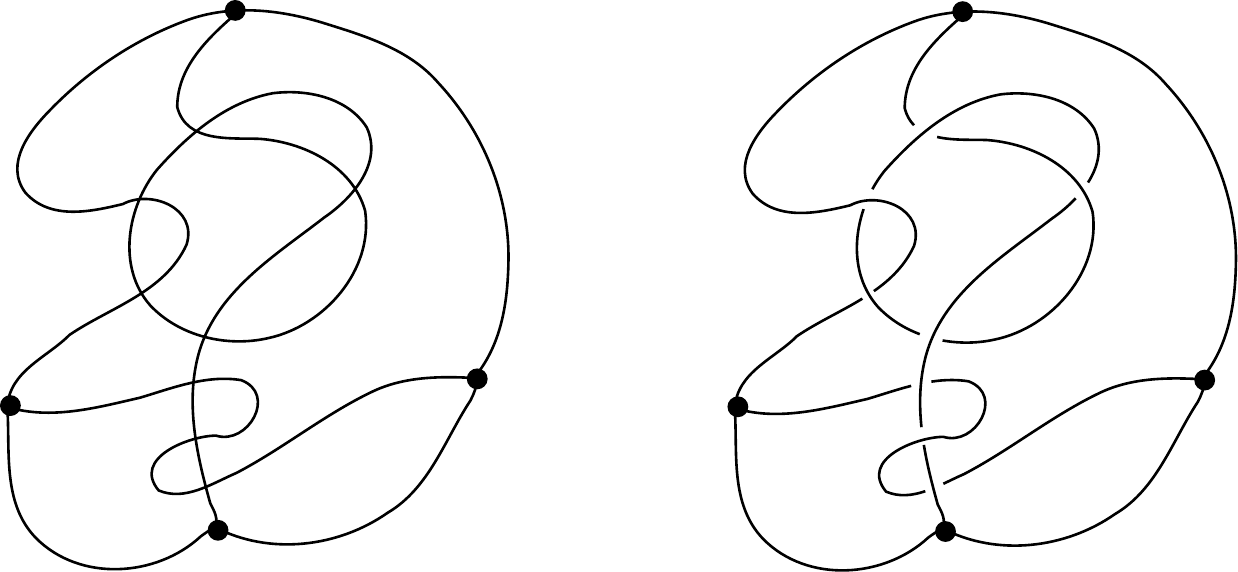}
\caption{Projection and Diagram}
\label{fig1}
\end{figure}

\begin{figure}
\centering
\includegraphics[width=9cm]{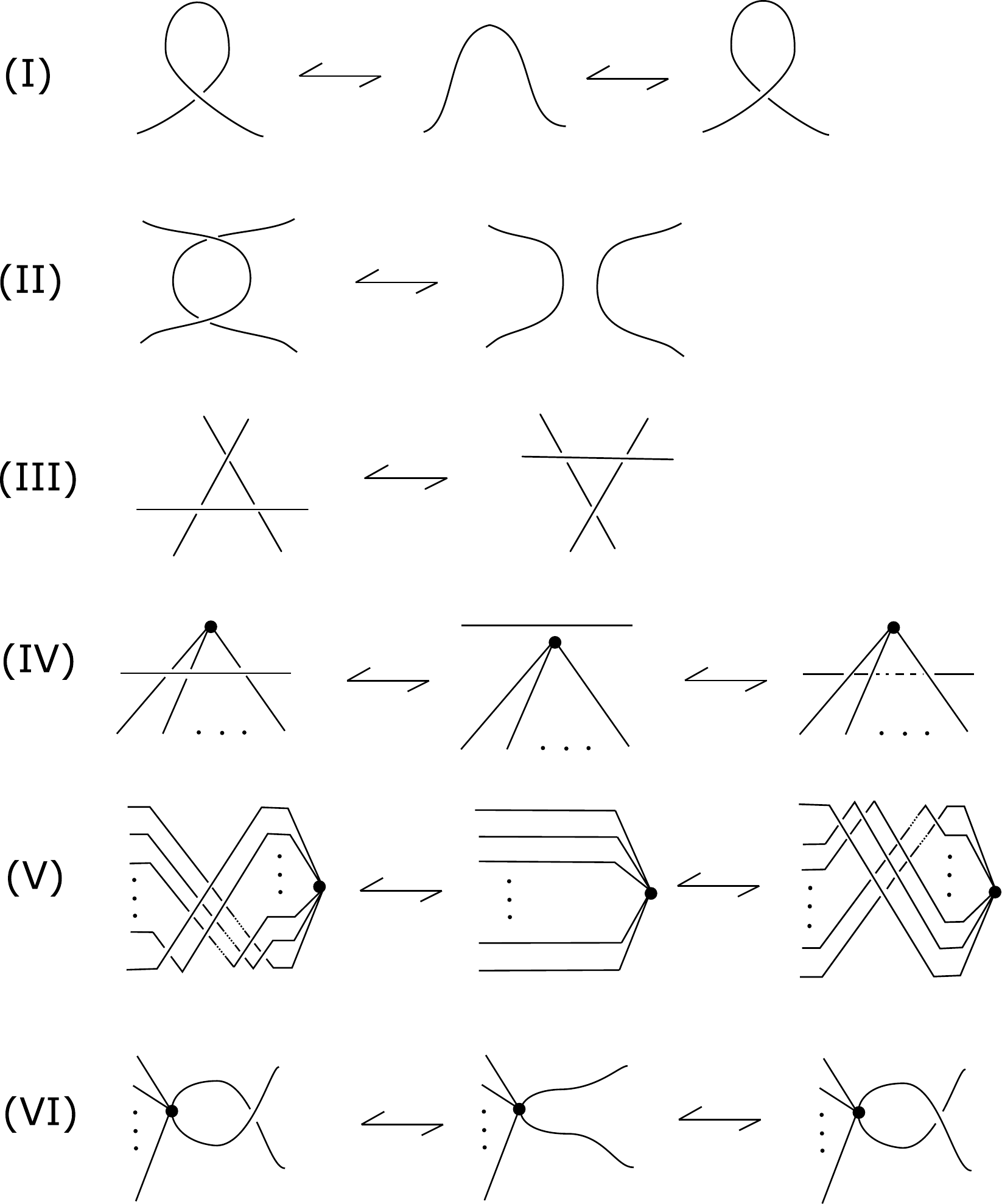}
\caption{Local moves on diagrams}
\label{fig2}
\end{figure}

Throughout this paper we work in the piecewise-linear category and graphs are considered to be 1-dimensional finite topological complexes. A {\em spatial graph} is a graph embedded into the Euclidean 3-space $\mathbb{R}^3$(or $\mathbb{S}^3$). Two spatial graphs $G$ and $G'$ are said to be {\em ambient isotopic}, if there exists an isotopy $h_t:\mathbb{R}^3 \rightarrow \mathbb{R}^3$ ($0\leq t\leq 1$) such that $h_0=id$ and $h_1(G)=G'$. A spatial graph $G$ is a {\em rigid-vertex graph}, if for each vertex $v$ there exists a neighborhood $N_v$ of $v$ in $G$ and a flat 2-disk $B_v$ centered at $v$ in $\mathbb{R}^3$ such that $G\cap B_v=N_v$. Two rigid-vertex graphs $G$ and $G'$ are said to be {\em rigid-vertex isotopic} (or shortly, {\em RV-isotopic}), if there exists an ambient isotopy $h_t$ such that $h_t(G)$ is also a flat-vertex graph for each $t$. These equivalences among spatial graphs can be described by the moves on diagrams. For a spatial graph $G$ let $p:\mathbb{R}^3 \rightarrow \mathbb{R}^2$ be a projection such that every multiple point of $p(G)$ is a transversal double point away from vertices. Then, adding information on over/underpassing at each crossing into the image $p(G)$, we obtain a {\em diagram} which represents the spatial graph $G$ (See Figure \ref{fig1}). Figure \ref{fig2} shows six local moves (I)$\sim$(VI) on diagrams.
\begin{lem}\cite{Ka} \label{lem1}
Two spatial graphs $G$ and $G'$ are ambient isotopic (resp. RV-isotopic) if and only if a diagram of $G$ can be transformed to a diagram of $G'$ by a finite sequence of moves among {\em (I)$\sim$(VI)} (resp. {\em (I)$\sim$(V)}).
\end{lem}
\begin{figure}
\centering
\includegraphics[width=9cm]{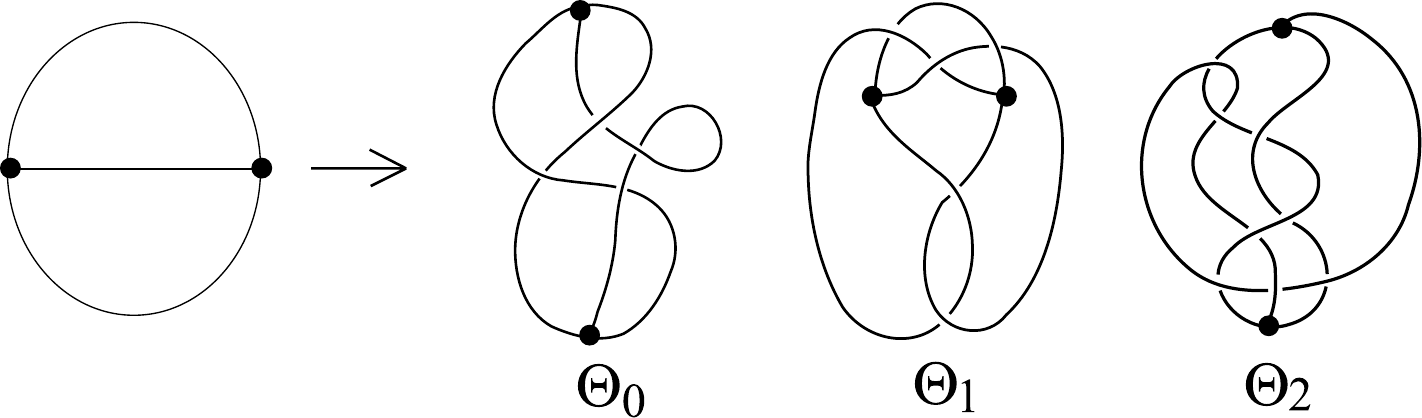}
\caption{$\theta$-graph and $\theta$-curves : $\Theta_0$ is trivial. $\Theta_1$ and $\Theta_2$ are brunnian.}
\label{fig3}
\end{figure}

In 1980's V. F. R. Jones discovered a polynomial invariant of knots and links \cite{Jones}. His discovery leads to the invention of more polynomial invariants of links. These invariants are strong and, furthermore, computable directly from diagrams of links. The idea and theory behind the invariants were generalized to the polynomial invariants of spatial graphs \cite{Jae, KV, Mu, Th, Ya, Yo}. The Yamada polynomial of spatial graphs \cite{Ya} is one of such invariants. The Yamada polynomial is a one-variable polynomial and a RV-isotopy invariant up to some mutiplication of its variable. Especially, for $\theta$-curves, the normalized Yamada polynomial is an ambient isotopy invariant and multiplicative under connected sum.

A {\em $\theta$-curve} is a spatial graph which consists of two vertices and three edges connecting the vertices (See Figure \ref{fig3}). A $\theta$-curve is said to be {\em trivial}\footnote{or {\em unknotted}}, if it is ambient isotopic to another $\theta$-curve which is contained in a plane. Note that a $\theta$-curve $\Theta$ contains three subknots $K_{\Theta}=\{e_i \cup e_j\;|\; i\neq j\}$, where $e_1$, $e_2$ and $e_3$ are the edges of $\Theta$. A non-trivial $\theta$-curve is said to be {\em brunnian}\footnote{Such $\theta$-curves are also said to be {\em almost trivial}, {\em locally unknotted} or {\em minimally knotted}.}, if every its subknot is of trivial knot type. Possibly knot invariants may distinguish two $\theta$-curves $\Theta$ and $\Theta'$ when $K_{\Theta}\neq K_{\Theta'}$ up to knot type. In the cases $K_{\Theta}= K_{\Theta'}$ we may use the Yamada polynomial. For example, the Yamada polynomial can distinguish $\Theta_1$ from $\Theta_0$ and $\Theta_2$ in Figure \ref{fig3}.

Links as well as knots can be utilized as invariants of spatial graphs. For a $\theta$-curve $\Theta$, let $S_{\Theta}$ be a closed two-punctured disk in $\mathbb{R}^3$ such that $\Theta$ is a spine of $S_{\Theta}$ and the Seifert form of $S_{\Theta}$ is zero. Then the 3-component link $L_{\Theta}$ corresponding to the boundary of $S_{\Theta}$ is called an {\em associated link} of $\Theta$.
Kauffman \textit{et al.} proved that the associated link of $\theta$-curves is ambient isotopy invariant of $\theta$-curves \cite{KSWZ}. The benefit of associated links is that we can use various link invariants as invariants of $\theta$-curves.
\begin{lem}\cite{KSWZ} \label{lem2}
\begin{enumerate}
\item For a $\theta$-curve $\Theta$, there exists a unique $S_{\Theta}$ up to ambient isotopy.
\item If two $\theta$-curves $\Theta$ and $\Theta'$ are ambient isotopic, then also $L_{\Theta}$ and $L_{\Theta'}$ are ambient isotopic.
\end{enumerate}
\end{lem}

In this paper we are interested in the relation between the Yamada polynomial of $\theta$-curves and the Jones polynomial of their associated links.
The Yamada polynomial, originally defined as a state sum on spatial graph diagrams \cite{Ya}, can be interpreted as a linear sum of the bracket polynomials of links obtained from some 2-parallels of spatial graph diagrams \cite{Jae}. Also the Jones polynomial of links can be reconstructed as the bracket polynomial \cite{Ka2}. Motivated by these facts, we investigate the difference between the normalized Yamada polynomial of $\theta$-curves and the normalized bracket polynomial of their associated links.

In \cite{Jae} Jaeger had defined a 2-variable polynomial of spatial graphs and proved that a specialization of his polynomial satisfies the skein relation of the Yamada polynomial. In Section 2 the Jaeger's version of the Yamada polynomial is briefly reviewed, and the precise relation between the original and the Jaeger's version is given.

In Section 3 we calculate the difference between the normaized Yamada polynomial of $\theta$-curves and the normalized bracket polynomial of associated links. As a corollary it is shown that the two polynomials coincide with each other for brunnian $\theta$-curves.

\section{Jaeger Polynomial and Yamada Polynomial}
In this section we briefly review the Jaeger polynomial and the Yamada polynomial (the readers are referred to \cite{Jae, Ya} for the details), and find the relation between them. Let $G$ denote an abstract graph, $\mathcal{G}$ an spatial embedding of $G$ (that is, spatial graph), and $D$ a diagram of $\mathcal{G}$.
\subsection{Jaeger Polynomial}
\begin{figure}
\centering
\includegraphics[width=9cm]{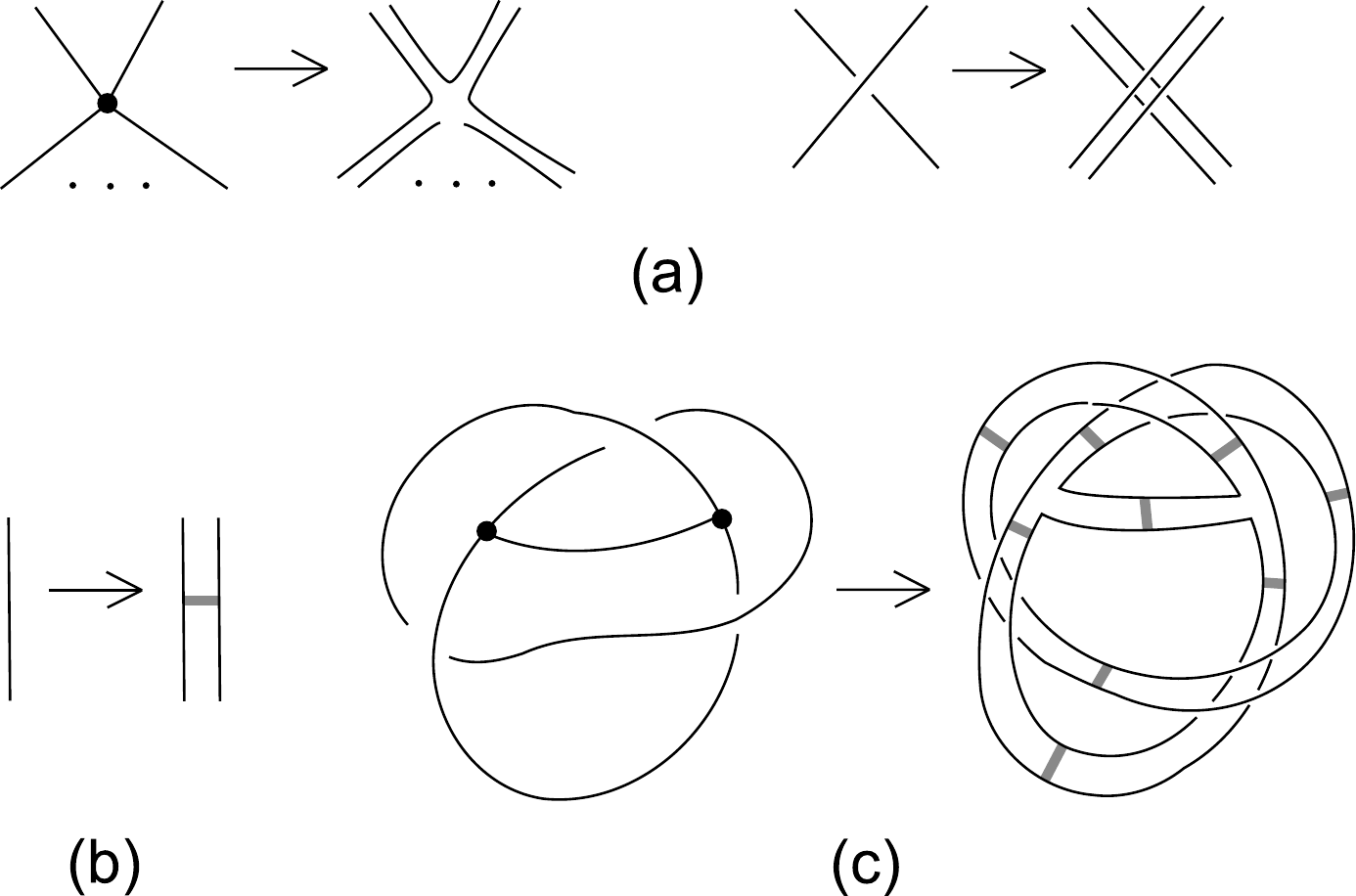}
\caption{Bar diagram}
\label{fig4}
\end{figure}

For a diagram $D$, the {\em bar diagram} $B_D$ is a band diagram with bars which is obtained from $D$. Figure \ref{fig4}-(a) shows the local parts of $B_D$ corresponding to  a vertex and a crossing of $D$, respectively. The other parts of $D$ (that is, edge strands of $D$ without vertices and crossings) are converted to strips with bars as illustrated in Figure \ref{fig4}-(b). Figure \ref{fig4}-(c) shows an example of bar diagram.
To the bar diagram $B_D$ we associate a fractional expression $R$ in two variables $a$ and $t$ which satisfies a defining relation:
$$R(\;\includegraphics[valign=c,width=0.4cm]{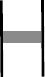}\;)
=\frac{1}{t+t^{-1}}\left(
R(\;\includegraphics[valign=c,width=0.4cm]{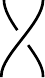}\;)
 + t^{-1}
R(\;\includegraphics[valign=c,width=0.4cm]{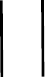}\;)
+ \frac{t-t^{-1}}{1-at}
R(\;\includegraphics[valign=c,width=0.4cm]{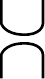}\;)
\right) $$
In the relation $R(B_D)$ is written in a combination of $R$'s of three bar diagrams. The three diagrams are obtained by changing one local strip with a bar of $B_D$ into a half-twist strip, a strip with no bar and a smoothing, respectively. Applying the relation repeatedly, $R(B_D)$ is written in terms of $R$'s of band diagrams with no bar. Then we can evaluate $R(B_D)$ by the second defining relation: {\em For a bar diagram $B$ with no bar,}
$$R(B)= \mathfrak{D}(L_B; a, z=t-t^{-1}), $$
where $L_B$ is a link diagram corresponding to the boundary of the band $B$, and $\mathfrak{D}(a,z)$ is the Dubrovnik polynomial of links by Kauffman\cite{Ka2}. For a diagram $D$ of a spatial graph $\mathcal{G}$, the Jaeger polynomial $J(D; a, t)$ is defined by
$$J(D)=R(B_D) \;.$$
Then the Jaeger polynomial is a RV-isotopic invariant upto powers of $at$.

The Jaeger's version of Yamada polynomial is a one-variable specialization of $J$:
$$\mathfrak{J}(D; A)= J(D, a=-A^3, t=A)\;.$$
Note that $\mathfrak{D}(a=-A^3,t=A)=\mathfrak{B}(A)$, where $\mathfrak{B}(L; A)$ is the bracket polynomial of link diagrams \cite{Ka2}.
\begin{lem}\cite{Jae} \label{lem3}
\begin{enumerate}
\item $\mathfrak{J}(A)$ is an RV-isotopy invariant of spatial graphs upto powers of $-A^4$.
\item Setting $a=-A^3$ and $t=A$,
$$R(\;\includegraphics[valign=c,width=0.4cm]{bar1}\;)
=
R(\;\includegraphics[valign=c,width=0.4cm]{bar3}\;)
+ \frac{1}{A^2+A^{-2}}
R(\;\includegraphics[valign=c,width=0.4cm]{bar4}\;) $$
\item For a spatial graph diagram $D$,
$$\mathfrak{J}(\;\includegraphics[valign=c,width=0.4cm]{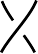}\;)
=A^4\; \mathfrak{J}(\;\includegraphics[valign=c,width=0.4cm]{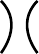}\;)
+ A^{-4}\; \mathfrak{J}(\;\includegraphics[valign=c,width=0.4cm]{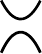}\;)
+ (A^2+A^{-2})\; \mathfrak{J}(\;\includegraphics[valign=c,width=0.4cm]{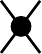}\;)
 $$
\end{enumerate}
\end{lem}
In the skein relation (3), the three diagrams in the right side are obtained by changing a crossing of $D$ into a vertical smoothing, a horizontal smoothing and a degree 4 vertex, respectively.

\subsection{Yamada Polynomial}
For an abstract graph $G$, let $V(G)$, $E(G)$, $\mu(G)$ and $\beta(G)$ be the set of vertices, the set of edges, the number of connected components and the first Betti number of $G$. Define a 2-variable polynomial $h$ of $G$ by
$$h(G; x, y)=\sum_{F\subset E(G)}(-x)^{-|F|}x^{\mu(G-F)}y^{\beta(G-F)} \;\;\; \mbox{and}\;\;\;  h(\emptyset)=1 \;.$$
For a spatial graph diagram $D$, a {\em state} is a function
$$s:\{ \mbox{crossings of } D \} \rightarrow \{+, -, 0\}\;.$$
Then let $D_s$ be a diagram obtained from $D$ by changing each crossing $c$ according to $s(c)$ as illustrated in the below. \\
\begin{center}
\includegraphics[width=5cm]{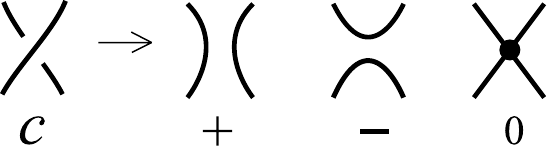}
\end{center}
Note that $D_s$ is a graph diagram with no crossing, hence represents a planar graph itself. Let $p(s)=|s^{-1}(+)|$ and $m(s)=|s^{-1}(-)|$. Then
the Yamada polynomial is defined by $Y(\emptyset)=1$ and
$$Y(D; A)= \sum_{s\in \mathcal{S}} A^{p(s)-m(s)}\;h(D_s; x=-1, y=-A-2-A^{-1}) \;, $$
where $\mathcal{S}$ is the set of states of $D$.
\begin{lem}\cite{Ya} \label{lem4}
\begin{enumerate}
\item $Y(A)$ is an RV-isotopy invariant of spatial graphs upto powers of $-A$.
\item For a spatial graph diagram $D$,
$$Y(\;\includegraphics[valign=c,width=0.4cm]{skein1}\;)
=A\; Y(\;\includegraphics[valign=c,width=0.4cm]{skein2}\;)
+ A^{-1}\;Y(\;\includegraphics[valign=c,width=0.4cm]{skein3}\;) +
Y(\;\includegraphics[valign=c,width=0.4cm]{skein4}\;)$$
\end{enumerate}
\end{lem}

\subsection{$\mathfrak{J}$ and $Y$}
The main purpose of this paper is to find the difference between the Yamada polynomial $Y$ of $\theta$-curves and the Jones polynomial of their associated links. The following proposition says that $\mathfrak{J}$ and $Y$ are equivalent, and hence allows us to observe $\mathfrak{J}$ for our purpose, instead of $Y$.
\begin{prop} \label{prop1}
Let $G$ be a planar graph and $D$ be a diagram of a spatial embedding of $G$. Then
$$Y(D; A^4)=-(A^2+A^{-2})^{|E(G)|-|V(G)|+1}\;\mathfrak{J}(D;A)$$
\end{prop}

\noindent {\bf Proof.} We prove the proposition by induction on the number of crossings of $D$. If $D$ has no crossing, then it represents the planar graph $G$ itself. It is pointed out by Jaeger \cite{Jae} that in this case $\mathfrak{J}$ is equivalent to the flow polynomial $F(G;t)$:
$$\mathfrak{J}(D;A)=(-A^2-A^{-2})^{|V(G)|-|E(G)|-1}F(G; A^4+2+A^{-4})\;.$$
The flow polynomial is a polynomial of abstract graphs and an invariant under graph isomorphism. It is a specialization of the Negami polynomial \cite{Ne} such that for a graph $H$
$$F(H;t)= t^{|E(H)|-|V(H)|}\sum_{F\subset E(H)} t^{\mu(H-F)}(-t)^{-|F|}\;.$$
Since $\beta=E-V+\mu$, we have
$$h(G; -1, y)=\sum_{F\subset E(G)}(-1)^{\mu(G-F)} y^{|E(G-F)|-|V(G-F)|+\mu(G-F)}\;.$$
Note that $|E(G-F)|=|E(G)|-|F|$ and $V(G-F)=V(G)$. Therefore
\begin{align*}
h(G; -1, y) &=  y^{|E(G)|-|V(G)|} \sum_{F\subset E(G)} (-y)^{\mu(G-F)} y^{-|F|}\\
&= (-1)^{|E(G)|-|V(G)|}(-y)^{|E(G)|-|V(G)|}\sum_{F\subset E(G)} (-y)^{\mu(G-F)} y^{-|F|}\\
&= (-1)^{|E(G)|-|V(G)|} F(G;-y)\;.
\end{align*}
Simply writing $E=E(G)$ and $V=V(G)$,
$$Y(D;\alpha)=h(G;-1,-\alpha-2-\alpha^{-1})
=(-1)^{|E|-|V|} F(G;\alpha+2+\alpha^{-1})\;.$$
Finally we have
\begin{align*}
Y(D;A^4) &= (-1)^{|E|-|V|} F(G;A^4+2+A^{-4})\\
&= (-1)^{|E|-|V|}(-A^2-A^{-2})^{|E|-|V|+1}\mathfrak{J}(G;A)\\
&=-(A^2+A^{-2})^{|E|-|V|+1}\mathfrak{J}(G;A) \;.
\end{align*}

Now suppose that $D$ has a crossing $c$. Let $D_{\infty}$, $D_{0}$ and $D_{*}$ be the diagrams obtained by changing the crossing
\includegraphics[valign=c,width=0.3cm]{skein1} into
\includegraphics[valign=c,width=0.3cm]{skein2},
\includegraphics[valign=c,width=0.3cm]{skein3} and
\includegraphics[valign=c,width=0.3cm]{skein4}, respectively.
Then the three diagrams have crossings one less than $D$.
In addition let $G_{\infty}$, $G_{0}$ and $G_{*}$ be the abstract graphs such that the three diagrams represent spatial embeddings of the three graphs, respectively. For our convenience write $\varepsilon=|E(G)|$ and $\nu=|V(G)|$.

Firstly consider the case that $c$ is a self-crossing, that is, a crossing produced by only one edge $e$. Then $G_{\infty}=G$,  and $G_{0}$ is the disjoint union of $G$ and one loop. The loop graph is considered to have one edge and one vertex. Therefore
$$|E(G_{\infty})|-|V(G_{\infty})|=\varepsilon-\nu=|E(G_{0})|-|V(G_{0})|\;.$$
The graph $G_{*}$ is the wedge sum of one loop and the subdivision on the edge $e$ of $G$, hence the difference is equal to $\varepsilon-\nu+1$. When $c$ is a crossing between two different edges, we have the same differences.

Applying the skein relation Lemma \ref{lem4}-(2) firstly, the induction hypothesis and then the difference enumeration,
\begin{align*}
Y(D;A^4) &= A^4\;Y(D_{\infty};A^4) + A^{-4}\; Y(D_0;A^4) + Y(D_*, A^4) \\
&= -A^4(A^2+A^{-2})^{\varepsilon-\nu+1} \mathfrak{J}(D_{\infty}; A) - A^{-4}(A^2+A^{-2})^{\varepsilon-\nu+1} \mathfrak{J}(D_{0}; A) \\
& \;\;\;\;\;\;\;\;\; - (A^2+A^{-2})^{\varepsilon-\nu+1+1} \mathfrak{J}(D_{*}; A)\\
&=  - (A^2+A^{-2})^{\varepsilon-\nu+1} \Big{\{} A^4\mathfrak{J}(D_{\infty}; A)+A^{-4}\mathfrak{J}(D_{0}; A)\\
& \;\;\;\;\;\;\;\;\; +(A^2+A^{-2})\mathfrak{J}(D_{*}; A) \Big{\}} \;.
\end{align*}
Finally, applying Lemma \ref{lem3}-(3), we complete the proof.
\qed
\section{Associated links of $\theta$-curves}
Let $D_{\Theta}$ be a diagram of a $\theta$-curve $\Theta$ with vertices $v_1$, $v_2$ and edges $e_1$, $e_2$, $e_3$. We may assume,  along the counter-clockwise direction, the edges appear in the order $(e_3, e_2, e_1)$ at $v_1$, and $(e_1, e_2, e_3)$ at $v_2$.
Thickening the edges in $D_{\Theta}$ we have a band diagram which represents a closed two punctured disk $S$. Then $\Theta$ is a spine of $S$.
Let $L=l_1 \cup l_2 \cup l_3$ be the 3-component link such that $\partial S=L$. Let the edges of $\Theta$ be oriented along the direction from $v_1$ to $v_2$. The link $L$ is oriented so that $l_1$ is homologous to $e_2-e_3$, $l_2$ to $e_3-e_1$ and $l_3$ to $e_1-e_2$ on $S$. The situation is depicted in Figure \ref{fig5}-(a).

\subsection{Associated links} To each crossing on $D_{\Theta}$ its sign is given by the following rule:
\begin{center}
\includegraphics[width=5cm]{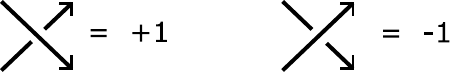}
\end{center}
\noindent Let $\omega_{ii}$ denote the sum of self-crossing signs of $e_i$, and $\omega_{ij}$ the sum of crossing signs between $e_i$ and $e_j$ ($i\neq j$). Then set three integers $n_i$ ($i=1, 2, 3$) to be
$$n_i= -\omega_{ii}+\frac{\omega_{ij}+\omega_{ik}-\omega_{jk}}{2}\;, \;\;\;\; \{i,j,k\}=\{1,2,3\}\;.$$

\begin{figure}
\centering
\includegraphics[width=9cm]{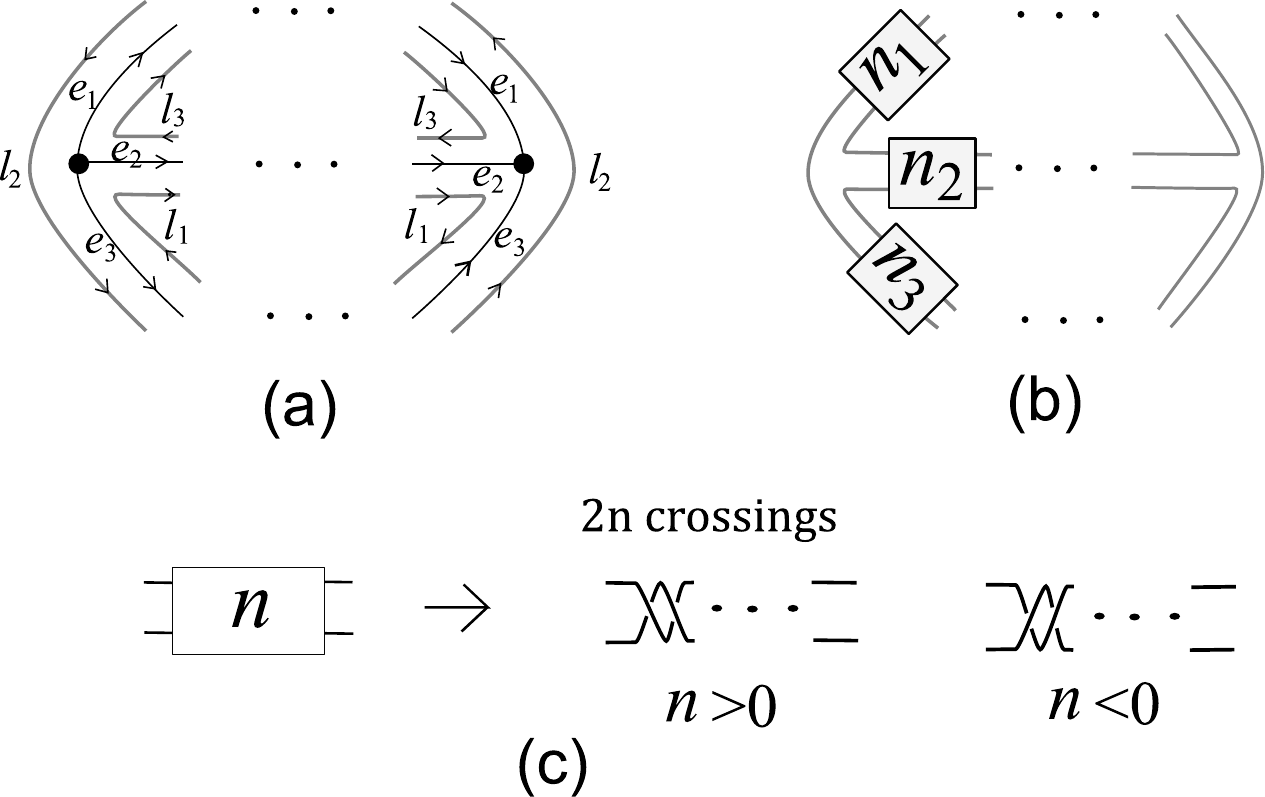}
\caption{Construction of associated link : (a) $\Theta$ and $L$, (b) associated link $L(n_1,n_2,n_3)$,  (c) n-full twist}
\label{fig5}
\end{figure}

Now, by adding three full-twist parts to $L$ as illustrated in Figure \ref{fig5}-(b) and (c), we obtain a link diagram $L(n_1,n_2,n_3)$ which represents the associated link of $\Theta$ \cite{KSWZ}. For our convenience let $<\;>$ denote the bracket polynomial with variable $A$, instead of $\mathfrak{B}(;A)$. Then the Jones polynomial $V(D)$ of an oriented link diagram $D$ corresponds to the normalization of $<D>$ \cite{Ka2}:
$$V(D;A)=(-A^3)^{-\omega}<D>\;,$$
where $\omega$ is the sum of all crossing signs on $D$. For our associated link $L(n_1, n_2, n_3)$, we know easily that
$$\omega=-2(n_1+n_2+n_3)\;.$$
\begin{prop} \label{prop2}
\begin{multline*}
V(L(n_1,n_2,n_3))= A^{8(n_1+n_2+n_3)}
\Big{\{} <L> + \sum_{i=1}^{3}\frac{1-(A^{-8})^{n_i}}{\varphi}<l_i^{(2)}> \\
+\frac{1}{\varphi^2}(2-\sum_{i=1}^{3}A^{-8n_i}+A^{-8(n_1+n_2+n_3)})
\Big{\}} \;,
\end{multline*}
where $\varphi=A^2+A^{-2}$ and $l_i^{(2)}$ is the 2-parallel of $l_{i}$.
\end{prop}

\noindent {\bf Proof.}
Firstly we give the properties of the bracket polynomial \cite{Ka2} which are necessary for our calculation:
\begin{itemize}
\item[(B0)] The bracket polynomial is invariant under moves (II) and (III).
\item[(B1)] $< \includegraphics[valign=c,width=0.4cm]{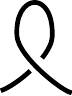}>=-A^3<\includegraphics[valign=c,width=0.4cm]{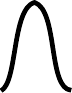}>$, $<\includegraphics[valign=c,width=0.4cm]{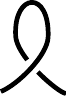}>=(-A)^{-3}<\includegraphics[valign=c,width=0.4cm]{curl-zero}>$, $<\includegraphics[valign=b,width=0.4cm]{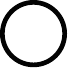}>=1$
\item[(B2)] $<D \cup \includegraphics[valign=b,width=0.4cm]{circ}>=(-A^2-A^{-2})<D>$, where $D\cup\includegraphics[valign=b,width=0.4cm]{circ}$ is a disjoint union of a link diagram $D$ and a simple circle.
\item[(B3)] $<\;\includegraphics[valign=c,width=0.4cm]{skein1}\;>
=A\; <\;\includegraphics[valign=c,width=0.4cm]{skein2}\;>
+ A^{-1}\;<\;\includegraphics[valign=c,width=0.4cm]{skein3}\;> $
\end{itemize}

\noindent Applying (B3) and (B1), we have
\begin{align*}
  <\includegraphics[valign=c,width=1.5cm]{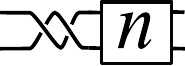}> &= A<\includegraphics[valign=c,width=1.25cm]{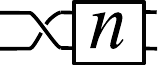}> + A^{-1}<\includegraphics[valign=c,width=1.3cm]{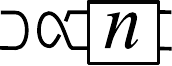}> \\
   &= A\Big{\{} A<\includegraphics[valign=c,width=0.77cm]{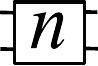}> + A^{-1}<\includegraphics[valign=c,width=1.2cm]{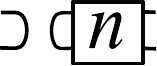}> \Big{\}} \\
   & \;\;\;\;\;\;\;\;\;\;\;\;\;\;\;\;\;\;\;\;\;\;\;\;\;\;\; +A^{-1}(-A^{-3})^{2n+1}<\includegraphics[valign=c,width=0.7cm]{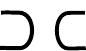}> \\
  &= A^{2}<\includegraphics[valign=c,width=0.77cm]{ftwist}> + \Big{\{}  (-A^{-3})^{2n} - A^{-4}(-A^{-3})^{2n} \Big{\}}<\includegraphics[valign=c,width=0.7cm]{ftwist-0s}>
\end{align*}
Assume $n>0$. Let $a_n=<\includegraphics[valign=c,width=0.77cm]{ftwist}>$ and $a_{\infty}=<\includegraphics[valign=c,width=0.7cm]{ftwist-0s}>$.
Then from the above equation
$$a_n=A^2 \;a_{n-1} + (1-A^{-4})(A^{-6})^{n-1} \;a_{\infty} \;.$$
Applying the reccurrence formula iteratively,
\begin{align*}
  a_n &= A^{2n}\;a_0 + (1-A^{-4})\sum_{k=1}^{n}A^{2(k-1)}(A^{-6})^{n-k}\; a_{\infty} \\
  &= A^{2n}\;a_0 +(1-A^{-4})\frac{A^{2n-2}(1-(A^{-8})^n)}{1-A^{-8}}\; a_{\infty}
\end{align*}
Therefore we have
$$a_n=A^{2n}\;a_0 + f_n\;  a_{\infty},\;\;\;\; f_n=\frac{A^{2(n-1)}(1-A^{-8n})}{1+A^{-4}} \;\;\cdots (*)$$
Similarly we can check that the equation $(*)$ is also true for $n<0$. In addition it holds that
$$A^{2n}f_m+A^{-6m}f_n=f_{n+m}=A^{2m}f_n+A^{-6n}f_m \;.$$
By averaging the left and the right we have
$$f_{n+m}=\frac{A^{2n}}{2}\Big{\{}1+A^{-8n}\Big{\}}f_m +\frac{A^{2m}}{2}\Big{\{}1+A^{-8m}\Big{\}}f_n\;. \;\;\cdots (**)$$

Now we calculate $<L(n_1, n_2, n_3)>$. Note the equivalences
\begin{align*}
&L(\infty, n_2,n_3)=l_1^{(2)}(n_2,n_3)=l_1^{(2)}(n_2+n_3),\\
&L(0, \infty, n_3)=l_2^{(2)}(n_3),\;\;
L(0,0,0)=L,\;\;
L(0,0,\infty)=l_3^{(2)}\;.
\end{align*}
Apply the equation $(*)$ to $n_1$, $n_2$ and $n_3$-full twists. Then
\begin{align*}
<L(n_1, n_2,n_3)>&=A^{2n_1}<L(0,n_2,n_3)> +f_{n_1}<L(\infty,n_2,n_3)>\\
&= A^{2n_1}\Big{\{} A^{2n_2}<L(0,0,n_3)>+f_{n_2}<L(0,\infty,n_3)>\Big{\}}\\
&\;\;\;\;\;\; +f_{n_1}<l_{1}^{(2)}(n_2+n_3)>\\
&= A^{2(n_1+n_2)}\Big{\{}A^{2n_3}<L(0,0,0)>+f_{n_3}<L(0,0,\infty)>\Big{\}} \\
&\;\;\;\;\;\; +A^{2n_1}f_{n_2}<l_2^{(2)}(n_3)> +f_{n_1}<l_1^{(2)}(n_2+n_3)>\\
&=A^{2(n_1+n_2+n_3)}<L>+A^{2(n_1+n_2)}f_{n_3}<l_3^{(2)}>\\
&\;\;\;\;\;\; +A^{2n_1}f_{n_2}<l_2^{(2)}(n_3)>+f_{n_1}<l_1^{(2)}(n_2+n_3)>\;.
\end{align*}
Again by $(*)$,
\begin{align*}
&<l_2^{(2)}(n_3)>=A^{2n_3}<l_2^{(2)}>+\;f_{n_3},\\
&<l_1^{(2)}(n_2+n_3)>=A^{2(n_2+n_3)}<l_1^{(2)}>+\;f_{n_2+n_3}\;.
\end{align*}
Therefore we have
\begin{multline*}
<L(n_1,n_2,n_3)>=A^{2(n_1+n_2+n_3)}\Big{\{} <L> + \sum_{i=1}^{3} \frac{f_{n_i}}{A^{2n_i}}<l_i^{(2)}> \\
+\frac{f_{n_2}}{A^{2n_2}}\frac{f_{n_3}}{A^{2n_3}} + \frac{f_{n_1}}{A^{2n_1}}\frac{f_{n_2+n_3}}{A^{2(n_2+n_3)}}
\Big{\}}\;,
\end{multline*}
\begin{align*}
V(L(n_1,n_2,n_3))&=(-A^3)^{2(n_1+n_2+n_3)}<L(n_1,n_2,n_3)>\\
&=A^{8(n_1+n_2+n_3)} \Big{\{}
<L> + \sum_{i=1}^3\frac{f_{n_i}}{A^{2n_i}}<l_i^{(2)}> \\
& \;\;\;\;\;\;\;\;+\frac{f_{n_2}}{A^{2n_2}}\frac{f_{n_3}}{A^{2n_3}} + \frac{f_{n_1}}{A^{2n_1}}\frac{f_{n_2+n_3}}{A^{2(n_2+n_3)}}
\Big{\}}\;.
\end{align*}
By the equation $(**)$,
\begin{align*}
&\frac{f_{n_1}}{A^{2n_1}}\frac{f_{n_2+n_3}}{A^{2n_2+2n_3}}\\
&=\frac{f_{n_1}}{A^{2n_1}}\frac{1}{A^{2n_2+2n_3}} \Big{\{}
\frac{A^{2n_2}}{2}(1+A^{-8n_2})f_{n_3} +
\frac{A^{2n_3}}{2}(1+A^{-8n_3})f_{n_2} \Big{\}} \\
&= \frac{1+A^{-8n_2}}{2} \frac{f_{n_1}}{A^{2n_1}} \frac{f_{n_3}}{A^{2n_3}} + \frac{1+A^{-8n_3}}{2} \frac{f_{n_1}}{A^{2n_1}} \frac{f_{n_2}}{A^{2n_2}}\;.
\end{align*}
Finally, the proof is completed by combining the above equation with $$\frac{f_n}{A^{2n}}=\frac{1-A^{-8n}}{\varphi}\;.$$
\begin{align*}
&\frac{f_{n_2}}{A^{2n_2}}\frac{f_{n_3}}{A^{2n_3}} + \frac{f_{n_1}}{A^{2n_1}}\frac{f_{n_2+n_3}}{A^{2(n_2+n_3)}}\\
&=\frac{1}{2\varphi^2}\Big{\{}
2(1-A^{-8n_2})(1-A^{-8n_3}) +(1+A^{-8n_2})(1-A^{-8n_1})(1-A^{-8n_3}) \\
&\;\;\;\;\;\;\;\;\;\;\;\;\;\;\;+(1+A^{-8n_3})(1-A^{-8n_1})(1-A^{-8n_2}) \Big{\}} \\
&=\frac{1}{2\varphi^2}\Big{\{} 4-2\sum_{i=1}^3 A^{-8n_i} + 2A^{-8(n_1+n_2+n_3)} \Big{\}}\;.
\end{align*}
\qed

\subsection{ $\mathfrak{J}$ of $\theta$-curves}
The normalized Yamada polynomial $$\widetilde{Y}(D_{\Theta})=(-A)^{2(n_1+n_2+n_3)}Y(D_{\Theta})$$ is an ambient isotopic invariant of $\theta$-curves \cite{Ya, HJ}. Hence, from Proposition \ref{prop1}, we directly know that the normalized Jaeger polynomial
$$\widetilde{\mathfrak{J}}(D_{\Theta})= (-A^4)^{2(n_1+n_2+n_3)}\mathfrak{J}(D_{\Theta})$$
is also an ambient isotopic invariant.
\begin{prop} \label{prop3}
$$\widetilde{\mathfrak{J}}(D_{\Theta})= (-A^4)^{2(n_1+n_2+n_3)}\Big{\{}
<L> + \frac{1}{\varphi}\sum_{i=1}^3 <l_i^{(2)}>+\frac{2}{\varphi^2} \Big{\}} $$
\end{prop}
\noindent {\bf Proof.} We calculate $\mathfrak{J}(D_{\Theta})=R(B_{D_{\Theta}};a=-A^3, t=A)$, where $B_{D_{\Theta}}$ is the bar diagram of $D_{\Theta}$. Firstly we give some properties of $R$ necessary for our calculation \cite{Jae}:
\begin{itemize}
\item[(R1)] $R( \;\includegraphics[valign=c,width=0.4cm]{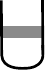}\;)=0$
\item[(R2)] $R(\; \includegraphics[valign=c,width=0.4cm]{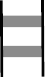}\;) = R(\includegraphics[valign=c,width=0.4cm]{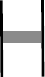}\;)$
\item[(R3)] $R(\;\includegraphics[valign=c,width=0.7cm]{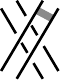}\;)=R(\;\includegraphics[valign=c,width=0.7cm]{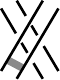}\;)$
\end{itemize}

\begin{figure}
\centering
\includegraphics[width=9cm]{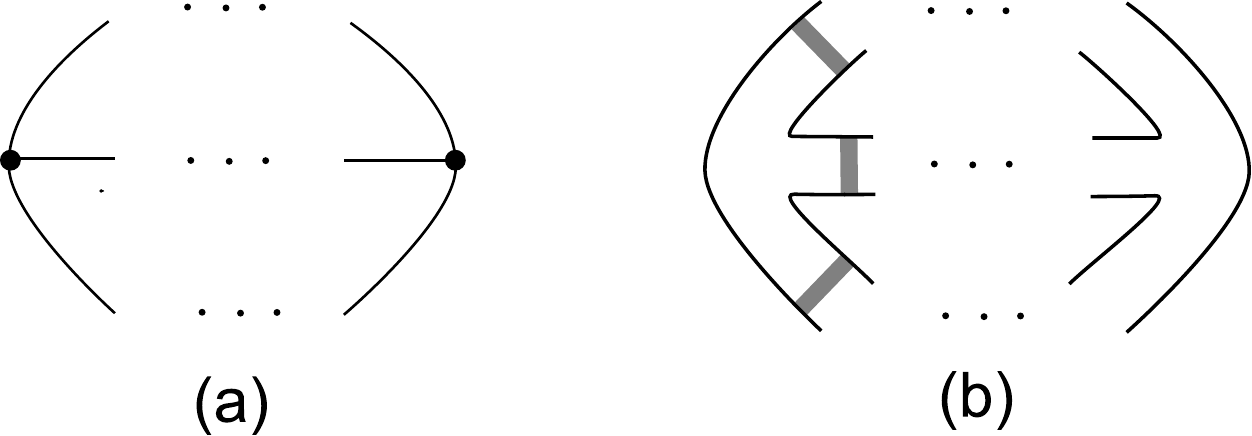}
\caption{(a) $D_{\Theta}$: a diagram of $\Theta$. (b) $B_{D_{\Theta}}$: $L$ with three bars.}
\label{fig6}
\end{figure}

For the $\theta$-curve diagram $D_{\Theta}$ the bar diagram $B_{D_{\Theta}}$ is simply the link diagram $L$ with bars. By the property (R3) it is allowed that the bars intersecting an edge $e_i$ are moved along the edge and positioned near the vertex $v_1$. Then the property (R2) allows us to merge the bars to be one bar.
In consequence $B_{D_{\Theta}}$ can be assumed to be $L$ with three bars, that is, only one bar for each edge, as depicted in Figure \ref{fig6}.
Let $L[\epsilon_1, \epsilon_2, \epsilon_3]$ $(\epsilon_i \in \{0, 1\})$ denote the link diagram $L$ with $\epsilon_i$ bar for the edge $e_i$. Therefore $L[0,0,0]=L$ and $L[1,1,1]=B_{D_{\Theta}}$. In addition let $l_i^{(2)}[n]$ denote the 2-parallel $l_i^{(2)}$ with $n$ bars, and let $U$ denote the simple circle $\includegraphics[valign=b,width=0.4cm]{circ}$.

To calculate $R(B_{D_{\Theta}})$, we apply the relation in Lemma \ref{lem3}-(2) iteratively. Then,
\begin{align*}
& R(L[1,1,1]) \\
&= R(L[0,1,1])+\varphi^{-1}R(l_1^{(2)}[2])\\
&= \Big{\{} R(L[0,0,1]) + \varphi^{-1}R(l_2^{(2)}[1]) \Big{\}} + \varphi^{-1}R(l_1^{(2)}[1])\\
&= \Big{\{} \Big( R(L[0,0,0])+\varphi^{-1}R(l_3^{(2)}) \Big)
+ \varphi^{-1}\Big(R(l_2^{(2)})+\varphi^{-1}R(U)\Big) \Big{\}}\\
& \;\;\;\;\;\; + \varphi^{-1}\Big{\{} R(l_1^{(2)})+\varphi^{-1}R[U] \Big{\}}\;.
\end{align*}
By definition, $$R(L[0,0,0])=R(L)=<L>,\; R(l_i^{(2)})=<l_i^{(2)}>,\; R(U)=<U>=1 \;.$$
Therefore we have
$$R(B_{D_{\Theta}})
= <L> + \frac{1}{\varphi}\sum_{i=1}^3 <l_i^{(2)}>+\frac{2}{\varphi^2} \;. $$
\qed

\section{Main Results}
\begin{thm} \label{thm1}
Let $\Theta$ be a $\theta$-curve with subknots  $\{\mathcal{K}_1, \mathcal{K}_2,\mathcal{K}_3\}$, and $\mathcal{L}$ be its associated link. Then,
$$\widetilde{\mathfrak{J}}(\Theta)-V(\mathcal{L})= \frac{1}{\varphi}\sum_{i=1}^3 \widetilde{\mathfrak{J}}(\mathcal{K}_i)\;-\frac{1}{\varphi^2}\;,$$
where $\varphi=A^2+A^{-2}$.
\end{thm}
From the theorem we know that the two polynomial invariants are actually equivalent each other for brunnian $\theta$-curves.
\begin{cor}\label{cor1}
If $\Theta$ is brunnian, then
$$\widetilde{\mathfrak{J}}(\Theta)-V(\mathcal{L})= -3 + \frac{2}{\varphi^2}\;.$$
\end{cor}

\vspace{0.2cm}
\noindent {\bf Proof of Theorem \ref{thm1}.}
Let $K$ be a diagram of a knot $\mathcal{K}$.
Then
\begin{align*}
\mathfrak{J}(K)&=R(K^{(2)}[1])\\
&=R(K^{(2)})+\frac{1}{\varphi}R(\includegraphics[valign=b,width=0.4cm]{circ})
\;=\;<K^{(2)}>+\frac{1}{\varphi}\;.
\end{align*}
The polynomial $\mathfrak{J}$ is invariant under moves (II) and (III), and
$$\mathfrak{J}(\includegraphics[valign=c,width=0.4cm]{curl-p})=
A^8\mathfrak{J}(\includegraphics[valign=c,width=0.4cm]{curl-zero})\;, \;\;\;\; \mathfrak{J}(\includegraphics[valign=c,width=0.4cm]{curl-m})=
A^{-8}\mathfrak{J}(\includegraphics[valign=c,width=0.4cm]{curl-zero})\;.$$
Therefore the normalized polynomial $\widetilde{\mathfrak{J}}$ is an ambient isotopy invariant of knots and we can write
$$\widetilde{\mathfrak{J}}(\mathcal{K})=\widetilde{\mathfrak{J}}(K)= (A^8)^{-\omega(K)}\mathfrak{J}(K)\;,$$
where $\omega(K)$ is the sum of crossing signs on $K$.

From Proposition \ref{prop2} and \ref{prop3},
\begin{align*}
&\widetilde{\mathfrak{J}}(\Theta)-V(\mathcal{L})\\
&=\widetilde{\mathfrak{J}}(D_{\Theta})-V(L(n_1,n_2,n_3))\\
&=A^{8(n_1+n_2+n_3)}\Big{\{}  \sum_{i=1}^{3}\frac{(A^{-8})^{n_i}}{\varphi}<l_i^{(2)}>
+\sum_{i=1}^{3}\frac{A^{-8n_i}}{\varphi^2}
\Big{\}}-\frac{1}{\varphi^2}\ \\
&=\frac{1}{\varphi}\Big{\{}
A^{8(n_2+n_3)}\Big(<l_1^{(2)}>+ \frac{1}{\varphi}\Big)
+A^{8(n_1+n_3)}\Big(<l_2^{(2)}>+ \frac{1}{\varphi}\Big)\\
&\;\;\;\;\;\; + A^{8(n_1+n_2)}\Big(<l_3^{(2)}>+ \frac{1}{\varphi}\Big)
\Big{\}} - \frac{1}{\varphi^2} \;.
\end{align*}
For the knot diagram $l_i$,
$$\omega(l_i)=\omega_{jj}+\omega_{kk}-\omega_{jk}=-(n_j+n_k) \;\;\;\;(\{i,j,k\}=\{1,2,3\})\;.$$
Therefore the difference is equal to
$$\frac{1}{\varphi}\sum_{i=1}^3 \widetilde{\mathfrak{J}}(l_i)\;-\frac{1}{\varphi^2}\;.$$
The knot represented by $l_i$ is of the same type with $\mathcal{K}_i$, hence $\widetilde{\mathfrak{J}}(l_i)=\widetilde{\mathfrak{J}}(\mathcal{K}_i)$.

\qed

\vspace{0.2cm}
\noindent {\bf Proof of Corollary \ref{cor1}.} Each subsknot $\mathcal{K}_i$ is trivial, hence
$$\widetilde{\mathfrak{J}}(\mathcal{K}_i) = \mathfrak{J}(\includegraphics[valign=b,width=0.4cm]{circ})=
-\varphi + \frac{1}{\varphi} \;.$$
Therefore the difference is equal to
$$\frac{1}{\varphi}\cdot 3 \cdot \Big( -\varphi + \frac{1}{\varphi}\Big) - \frac{1}{\varphi^2}=-3 + \frac{2}{\varphi^2}\; .$$
\qed
\section*{Acknowledgements}
The author was supported by the National Research Foundation of Korea (NRF) grant funded by the Korea government (MSIP).\\
(NRF - 2016R1D1A1B01008044)




\begin{thebibliography}{XX}
\bibitem{HJ} Y. Huh and G. T. Jin, {\em $\theta$-curve polynomials and finite-type invariants}, Knots 2000 Korea, Vol. 2 (Yongpyong). J. Knot Theory Ramifications \textbf{11} (2002), no. 4, 555–-564.
\bibitem{Jae} F. Jaeger, {\em On some graph invariants related to the Kauffman polynomial}, Progress in knot theory and related topics, 69–-82, Travaux en Cours, 56, Hermann, Paris, 1997.
\bibitem{Jones} V. F. R. Jones, {\em A polynomial invariant for knots via von Neumann algebras}, Bull. Amer. Math. Soc. (N.S.) \textbf{12} (1985), no. 1, 103–-111.

\bibitem{Ka} L. H. Kauffman, {\em Invariants of graphs in three-space}, Trans. Amer. Math. Soc. \textbf{311} (1989), no. 2, 697--710.

\bibitem{Ka2} L. H. Kauffman, {\em State models and the Jones polynomial}, Topology \textbf{26} (1987), no. 3, 395–-407.

\bibitem{KSWZ} L. H. Kauffman, J. Simon, K. Wolcott and P. Zhao, {\em
Invariants of theta-curves and other graphs in 3-space}, Topology Appl. \textbf{49} (1993), no. 3, 193–-216.

\bibitem{KV} L. H. Kauffman and P. Vogel, {\em Link polynomials and a graphical calculus}, J. Knot Theory Ramifications \textbf{1} (1992), no. 1, 59–-104.

\bibitem{Mu} J. Murakami, {\em The Yamada polynomial of spacial graphs and knit algebras}, Comm. Math. Phys. \textbf{155} (1993), no. 3, 511–-522.

\bibitem{Ne} S. Negami, {\em Polynomial invariants of graphs}, Trans. Amer. Math. Soc. \textbf{299} (1987), no. 2, 601–-622.

\bibitem{Th} A. Thompson, {\em A polynomial invariant of graphs in 3-manifolds}, Topology \textbf{31} (1992), no. 3, 657–-665.

\bibitem{Ya} S. Yamada, {\em An invariant of spatial graphs},
J. Graph Theory \textbf{13} (1989), no. 5, 537–551.

\bibitem{Yo} Y. Yokota, {\em Topological invariants of graphs in 3-space}, Topology \textbf{35} (1996), no. 1, 77–87.




\end{thebibliography}
\end{document}